\documentclass[graybox]{svmult}

\usepackage{helvet}         
\usepackage{courier}        
\usepackage{type1cm}        
\usepackage{makeidx}         
\usepackage{graphicx}        
\usepackage{multicol}        
\usepackage[bottom]{footmisc}

\usepackage{amsmath}
\usepackage{amssymb}
\usepackage{tikz}
\usepackage[enableskew]{youngtab}\Yboxdim{3pt} 

\newcommand{\tprod}[1][]{\textstyle\prod#1\displaystyle}
\newcommand{\tsum}[1][]{\textstyle\sum#1\displaystyle}

\newcommand{\C}{\mathbb{C}}
\newcommand{\bP}{\mathbb{P}}
\newcommand{\LG}[1]{\mathsf{#1}}
\newcommand{\LGD}{\LG{D}}

\newcommand{\OG}{\mathrm{OG}}
\newcommand{\SO}{\mathrm{SO}}
\newcommand{\Gr}{\mathrm{Gr}}
\newcommand{\PSO}{\mathrm{PSO}}

\newcommand{\G}{G}
\renewcommand{\P}{P}
\newcommand{\T}{T}

\newcommand{\cmX}{\mathbb{X}^\vee}
\newcommand{\mX}{X^\vee}
\newcommand{\omX}{X^\circ}
\newcommand{\dG}{G^\vee}
\newcommand{\dP}{P^\vee}
\newcommand{\dy}{y^\vee}

\newcommand{\can}{{\text{can}}}

\newcommand{\pot}{\mathcal{W}}

\newcommand{\de}{\delta}
\newcommand{\la}{\lambda}

\begin{document}

\title*{Canonical mirror models for maximal orthogonal Grassmannians}
\author{Peter Spacek and Charles Wang}
\institute{Peter Spacek \at TU Chemnitz, Germany, \email{peter.spacek@math.tu-chemnitz.de}
\and Charles Wang \at University of Michigan, Ann Arbor, US, \email{cmwa@umich.edu}}
\maketitle

\abstract{
In this report we use the methods of \cite{Spacek_Wang_Exceptional_family} and the upcoming generalization \cite{Spacek_Wang_Canonical_models} to complete the construction of canonical mirror models for all cominuscule homogeneous spaces, by considering the maximal orthogonal Grassmannians $\OG(n+1,2n+2)$.
}

\section{Definitions}\label{sec:definitions}
Let $X=\OG(n+1,2n+2)$ denote one of the two isomorphic connected components of the \emph{maximal orthogonal Grassmannian} of ($n+1$)-dimensional isotropic subspaces of $\C^{2n+2}$ with respect to the standard inner product. We consider $X$ as the homogeneous space $X=\G/\P$, where $\G=\SO(2n+2)$ is simple of type $\LGD_{n+1}$ and  $\P=\P_{n+1}$ is the maximal parabolic subgroup of $\G$ consisting of those $(2n+2)\times(2n+2)$-matrices that are block-upper-triangular when divided into four $(n+1)\times(n+1)$-blocks. This parabolic subgroup contains a natural choice of maximal torus $\T$, consisting of the diagonal elements of $\G$. Moreover, we will denote by $\widetilde{G}=\mathrm{Spin}(2n+2)$, $\widetilde{P}$ and $\widetilde{T}$ the universal covers of $\G$, $\P$ and $\T$, respectively.

We will define our \emph{canonical mirror models} as \emph{Landau-Ginzburg models} $(\mX_\can,\pot_\can)$, where $\mX_\can\subset \cmX$ is an open affine variety given by the complement of the singular locus of the \emph{superpotential} $\pot_\can: \cmX\times\C_q^*\dashrightarrow\C$ inside the \emph{Langlands dual} homogeneous space $\cmX=\dP\backslash\dG$. Here the \emph{Langlands dual groups} $\dG$ and $\dP$ of $\widetilde{G}$ and $\widetilde{P}$ are chosen relative to $\widetilde{T}$; in particular $\dG=\PSO(2n+2)$ and $\dP$ is the image of $\P$ under $\SO(2n+2)\to\PSO(2n+2)$. 

The superpotential $\pot_\can$ will be expressed in terms of \emph{(generalized) Pl\"ucker coordinates} obtained from the natural embedding $\cmX\hookrightarrow\bP(V_{n+1}^*)$ of $\cmX$ into the projectivization of the dual of the ($n+1$)st fundamental representation (also known as the \emph{spin representation}) of $\dG$.

The weight basis of $V_{n+1}^*$ induces \emph{Pl\"ucker coordinates} on $\bP(V_{n+1}^*)$ which can be labeled using (generalized) Young diagrams, analogously to the case for Grassmannians $\Gr(k,n)$.
A Young diagram%
\footnote{These diagrams arise naturally by considering \emph{order filters} on the \emph{minuscule poset} associated to the poset of weights of the fundamental representation $V_{n+1}$; such considerations work for any minuscule homogeneous space.}
for $V_{n+1}^*$ consists of a subset of the (``flipped'') staircase partition $\mu_n=(1,2,\dots, n)$, e.g. $\mu_3=\yng(1,2,3)$, such that each box has no ``empty spaces'' above or to the left of it inside of $\mu_n$. For $n=3$, the valid Young diagrams are $\varnothing$, $\raisebox{5pt}{\yng(1)}$, $\raisebox{2.5pt}{\yng(1,1)}$, $\raisebox{2.5pt}{\yng(1,2)}$, $\yng(1,1,1)$, $\yng(1,2,1)$, $\yng(1,2,2)$, $\yng(1,2,3)$, while $\yng(1,1,2)$ is invalid because there is an empty space above the second box in the last row. We also define $\la_i$ to be the diagram with $i$ columns of maximal length, e.g.~for $n=3$ we have $\la_1=\yng(1,1,1)$, $\la_2=\yng(1,2,2)$, $\la_3=\yng(1,2,3)$. Finally, we write $\la_i^+$ and $\mu_i^+$ for the diagrams obtained by adding the only possible box to $\la_i$ and $\mu_i$ respectively ($i<n$), e.g.~for $n=3$ we have $\la_1^+=\yng(1,2,1)=\mu_2^+$, and $\la_2^+=\mu_3$.

We will also define a rule to ``move'' boxes from one diagram to another. For this, we label the boxes with numbers from $1$ to $n+1$, starting with $1$ in the bottom left, $2$ in the boxes above and to the right of it, continuing to the last diagonal of boxes, which will be labeled alternatingly with $n+1$ and $n$ (starting with $n+1$ at the top). A box can only be \emph{removed} if it has no boxes to its right nor below it. On the other hand, a box labeled $i$ can only be \emph{added} if the diagram has an open spot labeled $i$ and the resulting diagram is still a valid Young diagram for $V_{n+1}^*$. 
We say that we can \emph{move} a box labeled $i$ from one diagram to another if a box with that label can be removed from one diagram and added to the other. To give a few examples of these operations, we consider $n=4$. In $\mu_2^+$, the boxes labeled $2$ and $4$ can be removed while boxes labeled $1$ and $3$ can be added; only the box labeled $4$ can be moved from $\mu_2^+$ to $\la_3$:
\vspace{-.5em}{
\begin{gather*}
\scriptsize\Yboxdim{7pt}
\raisebox{13.25pt}{\young(5,34)} ~~
\raisebox{14pt}{,}~~
\raisebox{6.5pt}{\young(5,3,2)}~~
\raisebox{14pt}{
$\xleftarrow[\textrm{2,4}]{\textrm{remove}}$
}~~
\raisebox{6.5pt}{\young(5,34,2)}~~
\raisebox{14pt}{$\xrightarrow[\textrm{1,3}]{\textrm{add}}$}~~
\young(5,34,2,1) ~~
\raisebox{14pt}{,}~~ 
\raisebox{6.5pt}{\young(5,34,23)}\quad
\raisebox{14pt}{;}\quad
\raisebox{6.5pt}{\young(5,34,2)} ~~ 
\young(5,34,235,123) ~~
\raisebox{14pt}{$\xrightarrow[4]{\textrm{move}}$} ~~
\raisebox{6.5pt}{\young(5,3,2)} ~~ \young(5,34,235,1234)
\end{gather*}\vspace{-2em}

}
\normalsize
\noindent As the labelling is always the same for a fixed $n$, we will suppress it in the following sections and write the diagrams only.

\section{The canonical mirror models}
The superpotential is written as $\pot_\can = \sum_{i=0}^{n+1}\pot_i$, where each $\pot_i$ is the quotient of homogeneous polynomials in Pl\"ucker coordinates of degree equal to the coefficient of the $i$th simple root in the highest root, except $\pot_0= \frac{p_{\yng(1)}}{p_{\!\varnothing}}$. For $i\in \{1,n,n+1\}$, the coefficient of the $i$th simple root in the highest root is equal to $1$, and the corresponding terms in our superpotential are:
\[
\pot_1 = \frac{p_{\la_1^+}}{p_{\la_1}}, \quad 
\pot_n = \frac{p_{\mu_{n-1}^+}}{p_{\mu_{n-1}}}, \quad\text{and}\quad 
\pot_{n+1} = q\frac{p_{\mu_{n-2}}}{p_{\mu_{n}}}.
\]
For $i\in [2, n-1]$, the coefficient of the $i$th simple root in the highest root is equal to $2$, and we will give an recursive definition for $\pot_i$ for $2\le i\le n-1$.

Beginning with $M_{i,0}=\{(\mu_{i-1},\la_i)\}$, we recursively define $M_{i,j}$ as the set of all pairs of Young diagrams that can be obtained from a pair $(\tau,\tau')\in M_{i,j-1}$ by moving a single box from $\tau$ to $\tau'$. Note that $\mu_{i-1}$ has $\ell_{i-1}=\binom{i}2$ boxes, so this process terminates. 

Next, define $N_{i,j}$ to be those pairs of Young diagrams that are obtained from pairs $(\tau,\tau')\in M_{i,j}$ by adding a box labeled $n+1-i$ to either $\tau$ or $\tau'$, if possible. It will follow from the type-independent results to be presented in \cite{Spacek_Wang_Canonical_models} that a box labeled $n+1-i$ can only be added to at most one of $\tau$ and $\tau'$. For $(\mu_{i-1},\la_i)\in M_{i,0}$ this can be seen directly: a box labeled $n+1-i$ can be added at the bottom of the first row of $\mu_{i-1}$, but $\la_i$ only accepts a box labeled $n$ or $n+1$ (depending on $i$ odd or even); hence, we always have $N_{i,0}=\{(\mu_{i-1}^+,\la_i)\}$. The term $\pot_i$ is then given by
\[
\pot_i = \frac
{\sum_{j}\sum_{(\tau,\tau')\in N_{i,j}} (-1)^j p_{\tau}p_{\tau'}}
{\sum_{j}\sum_{(\tau,\tau')\in M_{i,j}} (-1)^j p_{\tau}p_{\tau'}}.
\]
Finally, the canonical mirror variety is the complement $\mX_\can\subset\cmX$ of the singular locus of $\pot_\can$. This is the origin of the term \emph{canonical} mirror model: the degrees of the denominators of the terms of $\pot_\can$ add up to the index of $\cmX$ so they form an anti-canonical divisor and $\pot_\can$ is regular on the complement of an anti-canonical divisor.

\section{An overview of the proof}
We verify that $(\mX_\can,\pot_\can)$ provides a mirror model for $X$ by checking that each term $\pot_i$ agrees with a corresponding term of the Lie-theoretic potential defined in \cite{Rietsch_Mirror_Construction}. To do this, we show that the corresponding terms are equal after restriction to a certain algebraic torus $X^\circ$ which we now describe.

First, let $\dy_i:\C^*\to\dG$ be the one-parameter subgroup of $\dG$ associated to the negative of the $i$th simple root. Next, we set $\ell=\ell_n=\binom{n+1}{2}$ to be the number of boxes in $\mu_n$ and define $(i_1,\ldots,i_\ell)$ to be the sequence obtained by reading the labels of $\mu_n$ from left to right and top to bottom.%
\footnote{This can also be described as a choice of reduced expression for the minimal coset representative of $w_0W_\P$, where $W_\P$ is the Weyl group of (the Levi subgroup of) $\P$.} 
Furthermore, we also define $j_t$ to be the column that the label $i_t$ appears in. Now, we write $U_-^\circ\subset \dG$ for the open, dense torus of unipotent elements of the form $\dy_{i_\ell}(a_{i_\ell,j_\ell})\dots \dy_{i_{1}}(a_{i_1,j_1})$, and let $\omX = \dP\backslash U_-^\circ$.

We now compute the restrictions of the denominators\footnote{In \cite{Spacek_Wang_Exceptional_family}, we obtained these denominators from \emph{generalized minors} defined in \cite{GLS}, and we will show this holds generally in \cite{Spacek_Wang_Canonical_models}.} of each $\pot_i$, which we denote by $\phi_i$. We make use of the procedure presented in our previous paper \cite[Algorithm 3.12]{Spacek_Wang_Exceptional_family} to compute the restrictions of Pl\"ucker coordinates to the torus $\omX$. First, the restrictions of the denominators $\phi_0=p_{\varnothing}$, $\phi_1=p_{\la_1}$, $\phi_n=p_{\mu_{n-1}}$ and $\phi_{n+1}=p_{\mu_n}$ to $\omX$ are given by
\begin{align*}
\phi_0|_{\omX} &= 1, \qquad
\phi_1|_{\omX} = \tprod[_{t:j_t=1}] a_{i_t,j_t}, \\
\phi_{n}|_{\omX} &= \tprod[_{t=1}^{\ell_{n-1}}] a_{i_t,j_t},
\quad\text{and}\quad
\phi_{n+1}|_{\omX} = \tprod[_{t=1}^{\ell}] a_{i_t,j_t}.
\end{align*}
Note that $\phi_1|_{\omX}$ can be interpreted as the product of all $a_{i_t,j_t}$ with $i_t$ on the first column of $\mu_n$, while $\phi_n|_{\omX}$ can be interpreted as the product of all $a_{i_t,j_t}$ except for those with $i_t$ on the last row.
Next, we move to the denominators $\phi_i$ of $\pot_i$ for $i\notin\{1,n,n+1\}$. It will follow from the results in the upcoming \cite{Spacek_Wang_Canonical_models} that these are monomials in the $a_{i,j}$ when restricted to $\omX$, namely:
\[
\phi_i|_{\omX} = \left(\tprod[_{t=1}^{\ell_{i-1}}] a_{i_t,j_t}\right)\left(\tprod[_{t:j_t\le i}] a_{i_t,j_t}\right).
\]

We now compute the restrictions of the numerators of the $\pot_i$. We first consider $\pot_{n+1}$: comparing the steps of \cite[Algorithm 3.12]{Spacek_Wang_Exceptional_family} for $p_{\mu_{n-2}}$ restricted to $\omX$ with \cite[Corollary 8.12]{Spacek_Laurent_polynomial_LG_models}, we see that $\pot_{n+1}|_{\omX}$ is exactly the ``quantum term'' of the Laurent polynomial there (i.e.~the term with a factor of $q$). Next, we proceed to the numerators of $\pot_0$, $\pot_1$ and $\pot_n$: Note that each of these numerators is of the form $p_{\tau^{(i)}}$ where $\phi_i=p_{\tau}$ and $\tau^{(i)}$ is the tableau obtained from $\tau$ by adding a box labeled $n+1-i$ (which is a valid operation for the tableaux in question). Hence, we will write (for $i\in\{0,\ldots,n\}$)
\[
\delta_i(p_{\tau}) = 
\begin{cases}
    p_{\tau^{(i)}}, & \text{if $\tau^{(i)}$ is a valid Young diagram;}\\
    0, & \textrm{else,}
\end{cases}
\]
so that the numerator of $\pot_i$ (for $i\in\{0,1,n\}$) can be written as $\de_i(\phi_i)$.

Now, since the numerators are obtained from the denominators by adding a single box, labeled $n+1-i$, and the denominators are monomials on $\omX$, the crucial observation to obtain expressions for the numerators is that this box is the \emph{first} box with this label. By \cite[Algorithm 3.12]{Spacek_Wang_Exceptional_family} this has the implication that $\de_i(\phi_i)|_{\omX}$ is the product of $\phi_i|_{\omX}$ with $\sum_j a_{n+1-i,j}$, or in other words
\begin{equation}
\pot_i|_{\omX} = \tsum[_j] a_{n+1-i,j}.    
\label{eq:pot_i}
\end{equation}

For the numerators of $\pot_i$, where $i\notin\{0,1,n,n+1\}$, we extend the $\delta_i$ to monomials by the product rule $\delta_i(p_\tau p_{\tau'}) = \delta_i(p_\tau)p_{\tau'} + p_\tau\delta_i(p_{\tau'})$ and to polynomials by linearity. As before, the numerator of $\pot_i$ is exactly $\delta_i(\phi_i)$. To determine $\de_i(\phi_i)|_{\omX}$ for $i\notin\{0,1,n,n+1\}$, we use analogous arguments as for $i\in\{0,1,n\}$. However, one has to be more careful since the $\phi_i$ are now polynomials in Pl\"ucker coordinates. The same principle still holds due to the facts (i) that a monomial in $a_{i,j}$ in a given term $p_\tau p_{\tau'}$ of $\phi_i$ that cancels out against a monomial of another Pl\"ucker term will give rise to monomials in the numerator that cancel each other out as well, and (ii) that once $\de_i(p_\tau p_{\tau'})=0$ then $\de_i(p_{\tau''} p_{\tau'''})=0$ for all pairs $(\tau'',\tau''')$ obtained from $(\tau,\tau')$ by moving a series of boxes. These facts will be established in the upcoming \cite{Spacek_Wang_Canonical_models}.

Hence, we find that $\pot_i|_{\omX}$ satisfies \eqref{eq:pot_i} for 
all $i\in[0,n]$, implying that $\sum_{i=0}^n\pot_i|_{\omX} = \sum_{t=1}^\ell a_{i_t,j_t}$. Adding the term $\pot_{n+1}|_{\omX}$ on the left hand side completes the right hand side to the Laurent polynomial potential given in \cite[Corollary 8.12]{Spacek_Laurent_polynomial_LG_models}, implying in turn, by irreducibility of $\mX_\can$, that the canonical model is isomorphic to the Lie-theoretic model of \cite{Rietsch_Mirror_Construction}.

\section{An example}\label{sec:example}
For $n+1=4$, $\OG(4,8)$ is isomorphic to the six-dimensional quadric discussed in \cite{Pech_Rietsch_Williams_Quadrics}; one can verify that the potential defined there agrees with ours. Thus, we consider $n+1=5$, i.e.~$\OG(5,10)$. We draw the Hasse diagram of the poset of weights of the dual fundamental representation $V_5^*$ of $\dG$ labeled by the corresponding Young diagrams; the fully labeled $\mu_4$ is drawn to the right:
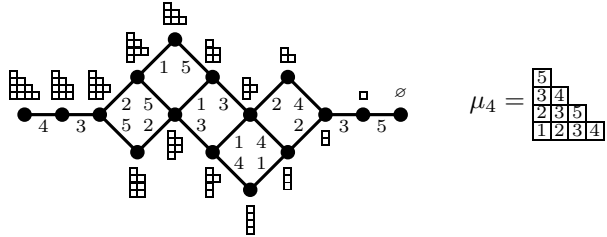
\vspace{-.5em}
\begin{equation*}
\begin{tikzpicture}[rotate=90,scale=.5,style=very thick,baseline=0.25em]
	\draw 
	(0,-7)--(0,-5)--(-2,-3)--(0,-1)--(-1,0)--(0,1)
    (0,-5)--(1,-4)--(0,-3)--(2,-1)
	(-1,-4)--(0,-3)--(-1,-2) 
	(0,1)--(2,-1) (1,-2)--(0,-1)--(1,0)
    (0,1)--(0,3)
    ;
	\draw[black,   fill=black] 
	(0,-7) circle (.15)
	(0,-6) circle (.15)
	(0,-5) circle (.15)
	(1,-4) circle (.15)
	(-1,-4) circle (.15)
	(0,-3) circle (.15)
	(-2,-3) circle (.15)
	(1,-2) circle (.15)
	(-1,-2) circle (.15)
	(2,-1) circle (.15)
	(0,-1) circle (.15)
	(1,0) circle (.15)
	(-1,0) circle (.15)
	(0,1) circle (.15)
    (0,2) circle (.15)
    (0,3) circle (.15)
    ;
	\node at (0,-6.5)[below=-2pt]{\scriptsize 5};
	\node at (0,-5.5)[below=-2pt]{\scriptsize 3};
	\node at (0.5,-4.5)[below left=-3pt]{\scriptsize 4};
	\node at (-0.5,-4.5)[above left=-3pt]{\scriptsize 2};
	\node at (0.5,-3.5)[below right=-3pt]{\scriptsize 2};
	\node at (-0.5,-3.5)[below left=-3pt]{\scriptsize 4};
	\node at (-1.5,-3.5)[above left=-3pt]{\scriptsize 1};
	\node at (0.5,-2.5)[below left=-3pt]{\scriptsize 3};
	\node at (-0.5,-2.5)[below right=-3pt]{\scriptsize 1};
	\node at (-1.5,-2.5)[above right=-3pt]{\scriptsize 4};
	\node at (1.5,-1.5)[below left=-3pt]{\scriptsize 5};
	\node at (0.5,-1.5)[below right=-3pt]{\scriptsize 1};
	\node at (-0.5,-1.5)[above right=-3pt]{\scriptsize 3};
	\node at (1.5,-0.5)[below right=-3pt]{\scriptsize 1};
	\node at (0.5,-0.5)[below left=-3pt]{\scriptsize 5};
	\node at (-0.5,0.5)[above right=-3pt]{\scriptsize 5};
    \node at (-0.5,-0.5)[above left=-3pt]{\scriptsize 2};
    \node at (0.5,0.5)[below right=-3pt]{\scriptsize 2};
    \node at (0,1.5)[below=-2pt]{\scriptsize 3};
    \node at (0,2.5)[below=-2pt]{\scriptsize 4};
	\node at (0,-7)[above=2pt]{\tiny$\varnothing$};
	\node at (0,-6)[above=2pt]{\tiny$\yng(1)$};
	\node at (0,-5)[below=2pt]{\tiny$\yng(1,1)$}; 
 \node at (1,-4)[above=2pt]{\tiny$\yng(1,2)$};
	\node at (-1,-4)[below=2pt]{\tiny$\yng(1,1,1)$};
	\node at (0,-3)[above=2pt]{\tiny$\yng(1,2,1)$};
	\node at (-2,-3)[below=2pt]{\tiny$\yng(1,1,1,1)$};
	\node at (1,-2)[above=2pt]{\tiny$\yng(1,2,2)$};
	\node at (-1,-2)[below=2pt]{\tiny$\yng(1,2,1,1)$};
	\node at (2,-1)[above=2pt]{\tiny$\yng(1,2,3)$};
	\node at (0,-1)[below=2pt]{\tiny$\yng(1,2,2,1)$};
	\node at (1,0)[above=2pt]{\tiny$\yng(1,2,3,1)$};
\node at (-1,0)[below=2pt]{\tiny$\yng(1,2,2,2)$};
\node at (0,1)[above=2pt]{\tiny$\yng(1,2,3,2)$};
\node at (0,2)[above=2pt]{\tiny$\yng(1,2,3,3)$};
\node at (0,3)[above=2pt]{\tiny$\yng(1,2,3,4)$};
\end{tikzpicture}
\qquad \mu_4=\scriptsize\Yboxdim{7pt}
\raisebox{-1.5em}{\young(5,34,235,1234)}
\normalsize
\vspace{-.5em}
\end{equation*}
Working out the definition of the superpotential, we obtain
\[
\pot_\can 
= \frac{p_{\raisebox{2.5pt}{\yng(1)}}\vphantom{p_{\yng(1,1,1,1)}}}{p_{\!\varnothing}}
+ \frac{p_{\yng(1,2,1,1)}}{p_{\yng(1,1,1,1)}}
+ \frac
    {p_{\yng(1,1)}p_{\yng(1,2,2,2)} - p_{\!\varnothing}p_{\yng(1,2,3,3)}}
    {p_{\raisebox{2.5pt}{\yng(1)}}p_{\yng(1,2,2,2)} - p_{\!\varnothing}p_{\yng(1,2,3,2)}}
+ \frac
    {p_{\yng(1,2,1)}p_{\yng(1,2,3,3)} - p_{\yng(1,1,1)}p_{\yng(1,2,3,4)}}
    {p_{\yng(1,2)}p_{\yng(1,2,3,3)}-p_{\yng(1,1)}p_{\yng(1,2,3,4)}}
+ \frac{p_{\yng(1,2,3,1)}}{p_{\yng(1,2,3)}}
+q\frac{p_{\yng(1,2)}\vphantom{p_{\yng(1,1,1,1)}}}{p_{\yng(1,2,3,4)}}
\]
and we single out the fourth term $\pot_3$ on the right-hand side for further consideration. First, we compute the restrictions of the relevant Pl\"ucker coordinates in the denominator $\phi_3$ of $\pot_3$ to $\omX$. For this, recall the labeling for the torus coordinates $a_{{i_t},{j_t}}$ given by 
\[
(i_t,j_t)_{t=1}^{10} = \bigl((5,1),(3,1),(4,2),(2,1),(3,2),(5,3),(1,1),(2,2),(3,3),(4,4)\bigr).
\] 
With this, \cite[Algorithm 3.12]{Spacek_Wang_Exceptional_family} tells us that the restrictions are given by:
\begin{align*}
p_{\yng(1,2)}|_{\omX} &= a_{5,1}(a_{3,1}a_{4,2}+a_{3,1}a_{4,4}+a_{3,2}a_{4,4}+a_{3,3}a_{4,4})+a_{5,3}a_{3,3}a_{4,4},\\
p_{\yng(1,2,3,3)}|_{\omX} &= a_{5,1}a_{3,1}a_{4,2}a_{2,1}a_{3,2}a_{5,3}a_{1,1}a_{2,2}a_{3,3},\\
p_{\yng(1,1)}|_{\omX} &= a_{5,1}(a_{3,1}+a_{3,2}+a_{3,3})+a_{5,3}a_{3,3},\\
p_{\yng(1,2,3,4)}|_{\omX} &= a_{5,1}a_{3,1}a_{4,2}a_{2,1}a_{3,2}a_{5,3}a_{1,1}a_{2,2}a_{3,3}a_{4,4}.
\end{align*}
Using the expressions computed above, we compute the restriction
\begin{equation}
\phi_3|_{\omX} = a_{5,1}^2a_{3,1}^2a_{4,2}^2a_{2,1}a_{3,2}a_{5,3}a_{1,1}a_{2,2}a_{3,3}.
\label{eq:D5_phi3_torus}
\end{equation}

Next, we verify that the numerator of $\pot_3$ is given by $\delta_3(\phi_3)$:
\[
\phi_3 = p_{\yng(1,2)}p_{\yng(1,2,3,3)}-p_{\yng(1,1)}p_{\yng(1,2,3,4)}
\quad\text{with labeling}\quad
\scriptsize\Yboxdim{7pt}
\raisebox{1.5pt}{\young(5,34)}\quad 
\raisebox{-12pt}{\young(5,34,235,123)}\quad 
\raisebox{1.5pt}{\young(5,3)}\quad 
\raisebox{-12pt}{\young(5,34,235,1234)}
\normalsize .
\]
A box labeled $4+1-3=2$ can only be added to $\raisebox{2.5pt}{\yng(1,1)}$ and $\raisebox{2.5pt}{\yng(1,2)}$, yielding $\yng(1,1,1)$ and $\yng(1,2,1)$ respectively. This means that $\de_3$ only acts nontrivially on $p_{\yng(1,1)}$ and $p_{\yng(1,2)}$. In particular, we have $\delta_3(\phi_3) = p_{\yng(1,2,1)}p_{\yng(1,2,3,3)} - p_{\yng(1,1,1)}p_{\yng(1,2,3,4)}$, which is indeed the numerator. Moreover,\vspace{-.5em}
\begin{align*}
p_{\yng(1,2,1)}|_{\omX} &= a_{5,1}\bigl(a_{3,1}a_{4,2}(a_{2,1}+a_{2,2})+a_{3,1}a_{4,4}(a_{2,1}+a_{2,2}) + a_{3,2}a_{4,4}a_{2,2}\bigr),\\
p_{\yng(1,1,1)}|_{\omX} &= a_{5,1}a_{3,1}(a_{2,1}+a_{2,2})+a_{5,1}a_{3,2}a_{2,2},
\end{align*}
so we compute the torus restriction of the numerator as
\[
\de_3(\phi_3)|_{\omX}
= a_{5,1}^{2}a_{3,1}^{2}a_{4,2}^{2}a_{2,1}a_{3,2}a_{5,3}a_{1,1}a_{2,2}a_{3,3}(a_{2,1}+a_{2,2})
\]
Note that this is exactly \eqref{eq:D5_phi3_torus} times $a_{2,1}+a_{2,2}$, so $\pot_3|_{\omX}=a_{2,1}+a_{2,2}$ as claimed: $a_{2,1}$ and $a_{2,2}$ correspond exactly to the two occurrences of $2$ in the sequence $(i_t)=(5,3,4,2,3,5,1,2,3,4)$ obtained from $\mu_4$.

\end{document}